\documentclass[letterpaper, 10 pt, conference]{article}
  \usepackage{amsfonts}

\usepackage{graphics}
\usepackage{graphicx}
\usepackage{epsfig} 
\usepackage{mathptmx} 
\usepackage{times} 
\usepackage{amsmath} 
\usepackage{amssymb}  
\usepackage{verbatim} 
\usepackage[linesnumbered,ruled,vlined]{algorithm2e}
\usepackage{subfigure}
\usepackage{url,cite}
\usepackage[linesnumbered,ruled,vlined]{algorithm2e}
\usepackage{color}

\newtheorem{definition}{Definition}
\newtheorem{remark}{Remark}

\newtheorem{assumption}{Assumption}

\newcommand{\Pcal}{\mathcal{P}}
\newcommand{\Hcal}{\mathcal{H}}

\newcommand{\jN}{j_{|N_i|}}
\newcommand{\zeros}{\textbf{0}}
\newcommand{\Id}{\textbf{I}}

\title{\LARGE \bf
Information-driven Fully Distributed Kalman Filter for Sensor Networks in Presence of Naive Nodes
}

\author{Shaocheng Wang, Wei Ren and Zhongkui Li
\thanks{Shaocheng Wang and Wei Ren are with the University of California, Riverside, CA 92507, USA. Email: shaocheng.wang@email.ucr.edu, ren@ee.ucr.edu.}
\thanks{Zhongkui Li is with the Peking University, P. R. China. Email: zhongkli@pku.edu.cn}
}

\begin{document}

\maketitle
\thispagestyle{empty}
\pagestyle{empty}

\begin{abstract}
We consider the distributed Kalman filtering problem for sensor networks where each node takes the measurement, communicates with its local neighbors, and updates its local estimate and estimation error covariance at the same frequency.
In such a scenario, if the target is not directly observed by neither a certain node nor its local neighbors, this node is naive about the target.
The well-known Kalman Consensus filter (KCF) has been shown to perform well if there exists no naive node in the network.
The case in presence of naive nodes has been considered by the generalized KCF (GKCF) and the Information Weighted Consensus filter (ICF) later on.
However, all these consensus-embedded filters require some global information such as the maximum degree of the graph, or the total number of the nodes.
With communication topology changes, node failures, or addition of new nodes, the filter performance would be adversely affected.
In this paper, we consider a novel local Weighted Least Square estimator for each node that utilizes its \emph{generalized measurement} formed by not only its own and local neighbors' measurements but their prior local estimates to track the target.
With some approximations in the derivation of the covariance matrix, we propose the Information-driven Fully Distributed Kalman filter (IFDKF),
which is able to deal with the existence of naive nodes without knowing any global information.
Experimental results show that the proposed algorithm performs better than the existing algorithms in the considered realistic scenario.

\end{abstract}

\section{Introduction}
Estimating and tracking the state of a dynamic process is one of the most fundamental but significant problems in sensor networks.
The Kalman filter\cite{Kalman1960} is known as the optimal solution, in the sense of the mean squared error (MSE), to the linear discrete filtering problems.
Traditional approaches of tracking a dynamic target usually require a centralized server, which collects the local measurements from each sensor at each time instant, and estimates the state of the target.
These centralized solutions usually guarantee the optimality of the estimates.
However, they also require a large number of communications, which could bring issues especially with the increase of the network size and the limited communication bandwidths.

Therefore, in recent years, there is an increasing interest in the distributed estimation and tracking problem \cite{RenBeardKingston05
,AlighanbariUnbiasedKalman,olfati2007distributedKalman,olfati2009kalman}.
Unlike the traditional approaches, the distributed estimation mechanisms try to recover the centralized solution via peer-to-peer communications.
That is, each node of the network tries to estimate the state of the target, by communicating with only its local neighbors.
This would not only save the communication cost in the network, but also improve the network robustness against the possible failure of sensors.
The consensus algorithm, which computes the global (weighted) average of a variable of interest in a distributed manner, was considered from the Kalman filters' perspective in \cite{RenBeardKingston05}.
Ref. \cite{AlighanbariUnbiasedKalman} proposed a modified version of the algorithm in \cite{RenBeardKingston05} and showed that the modified algorithm was unbiased even if the outflows of each node are not equal to each other.
Both algorithms in \cite{RenBeardKingston05} and \cite{AlighanbariUnbiasedKalman} were proposed to estimate a static state.
The well-known Kalman Consensus filter (KCF), as a combination of the consensus filter and the distributed Kalman filter, was proposed in \cite{olfati2007distributedKalman} to estimate a dynamic state.
At each node, the KCF sends/receives the measurements to/from its local neighbors and an average-consensus term is added in the measurement updates to guarantee that the estimates of each sensor asymptotically become consistent.
The optimality and stability of the KCF was analyzed in \cite{olfati2009kalman}.
Ref. \cite{bai2011distributed} applied the results of the dynamic averaging algorithm, a generalization of the average consensus algorithm, into the distributed Kalman filtering problem.

In the sensor networks, it is possible that a target is observed by neither a certain node nor any of its local neighbors.
This type of nodes was referred as the \emph{naive} nodes with respect to the target in \cite{kamal2011GKCF}.
The naive nodes usually have bad local estimates, which will be further sent to its neighbors and consequently contaminate the estimation of the whole network through the aforementioned average-consensus algorithm.
Therefore, the KCF performs well only if there exists no naive node in the network.
However, the network with no naive node is not always the case in real applications.
For example, in a network with multiple cameras designed to monitor a large area, as studied in \cite{kamal2011GKCF,kamal2013ICFTAC},
each camera tries to compute the state of a moving target while some cameras are naive about the target due to their limited field of view and local communication.
In such a scenario, the KCF will not be able to track the target as each node weighs its neighbors' estimates in an equal manner.
This issue was considered in \cite{kamal2011GKCF}, where the generalized Kalman Consensus filter (GKCF) was proposed.
The GKCF allows each node to weigh its neighbors' estimates by the inverse of their estimation error covariance matrices, and therefore, avoids the estimates from diverging in presence of naive nodes.

The information filter, which is mathematically equivalent to the Kalman filter, was combined with the consensus filter in \cite{casbeer2009distributed}.
It was shown in \cite{casbeer2009distributed} that the proposed filter made each node overestimate its estimation uncertainty by a factor of the number of nodes in the network.
This would degrade the performance with the increase of the network size.
Ref. \cite{kamal2013ICFTAC} proposed the Information Weighted Consensus filter (ICF), analyzed its properties, and showed the comparisons of its performance with the KCF and the GKCF.
The ICF was shown to approach the centralized optimal solution as the nodes communicate with their local neighbors for infinite times before they update their local estimates.
Both the GKCF and the ICF outperform the KCF in the presence of naive nodes.
However, both algorithms require some global information such as the maximum degree of the graph, the total number of nodes in the network, and/or proper selection of certain parameters (the step sizes in the consensus update).
This would be a limitation if such global information changes with time, or the step sizes in the consensus update are not properly selected.

In this paper, we propose an Information-driven Fully Distributed Kalman filter (IFDKF).
In the proposed algorithm, each node obtains its local measurement, communicates with its local neighbors, and updates its local estimate and estimation error covariance at the same frequency.
Different from the KCF, the IFDKF is able to work in presence of naive nodes.
More importantly, the IFDKF is fully distributed and therefore robust against the possible change of the aforementioned global parameters.
Each node runs it algorithm using only local information in an automated manner and there is no need to tune or select any parameter across the nodes in the algorithm.
The experimental results show that the proposed IFDKF outperforms other distributed Kalman-based filters.

\section{Preliminary}
\subsection{Graph Theory}
A graph $\mathcal{G(V,E)}$ is used to represent the topology of a sensor network, where $\mathcal{V}$ is the set of vertices that stands for the  sensors, and $\mathcal{E\subseteq V\times V}$ is the set of edges that stands for the communication channels, respectively.
Suppose that there are $N$ nodes and $l$ communication channels in the network.
We let $\mathcal{V}=\{v_1,...,v_{N}\}$.
Several basic concepts of graph theory used later are briefly listed here.
An \emph{edge} $(i,j)\in\mathcal{E}$ denotes that node $j$ can obtain information from node $i$.
If a graph is undirected, $(i,j)\in\mathcal{E}$ implies that $(j,i)\in\mathcal{E}$.
A \emph{path} from vertex $v_{i_0}$ to vertex $v_{i_\ell}$ is a sequence of vertices $v_{i_0},v_{i_1},...,v_{i_t}$ such that $({i_{j-1}},{i_{j}})\in\mathcal{E}$ for $0<j\leq \ell$.
The \emph{distance} from vertex $v_{i}$ to vertex $v_{j}$ is the number of sequences contained in the shortest path between them.
A graph is \emph{connected} if there exists at least one path from every vertex to every other vertex.
$N_i\triangleq\{j|(j,i)\in\mathcal{E},\forall j\neq i\}$ is the \emph{neighborhood} of node $i$.
$J_i\triangleq N_i\cup\{i\}$ is the \emph{inclusive neighborhood} of node $i$.
The \emph{degree} of a certain node is defined as $\Delta_i \triangleq |N_i|$.
The \emph{maximum degree} of a graph is defined as $\Delta_\textrm{max} \triangleq \max_{i\in\{1,\ldots,N\}} \Delta_i$.
A \emph{spanning tree}\footnote{Here we focus on the spanning tree in an undirected graph.} is a minimal set of edges that connect all vertices.

\begin{figure}
\centering
\subfigure[]{\label{fig:UndirectGraph1} 
\includegraphics[width=.3\linewidth]{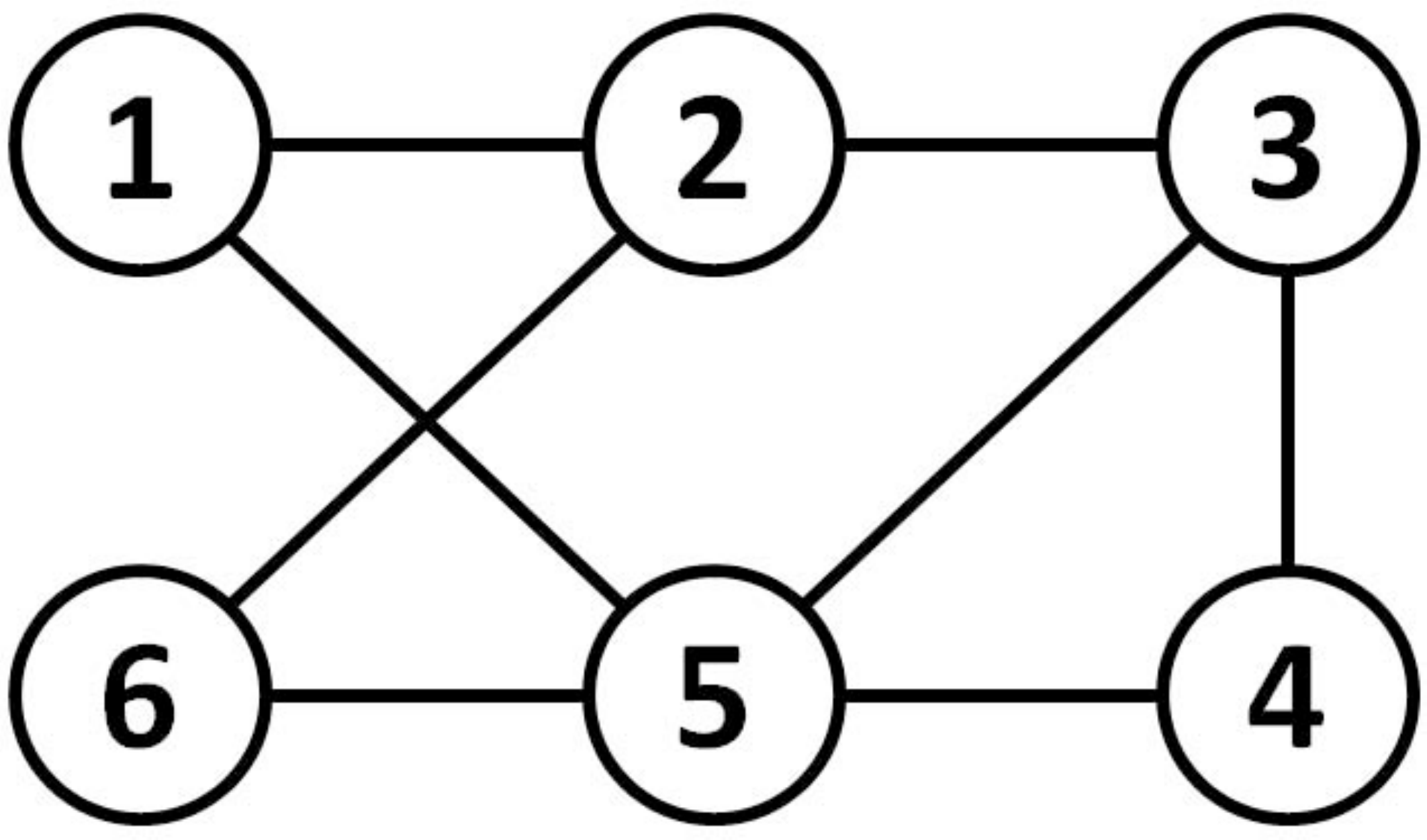}
}
\subfigure[]{\label{fig:UndirectGraph2} 
\includegraphics[width=.3\linewidth]{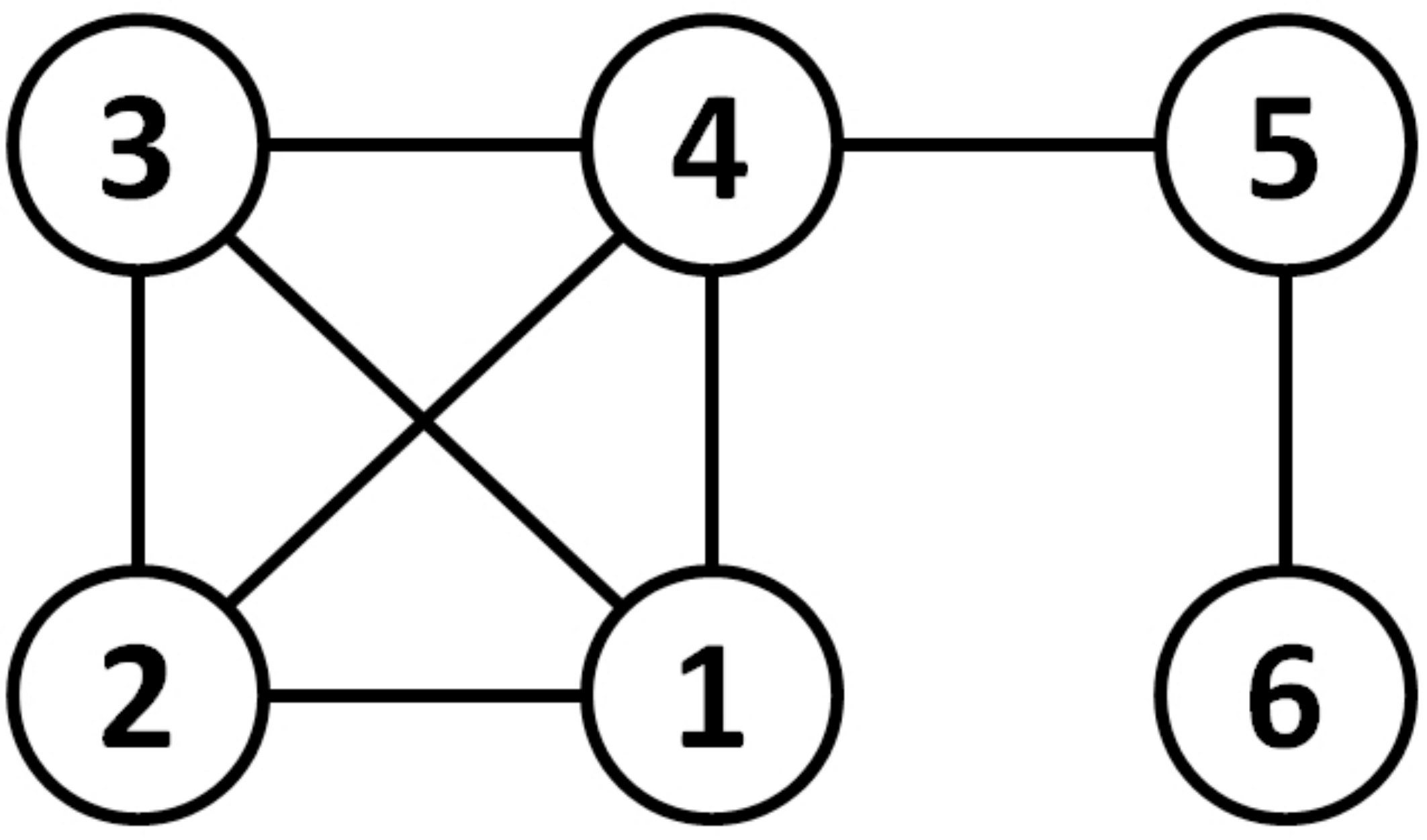}}
\subfigure[]{\label{fig:SpanningTree} 
\includegraphics[width=.3\linewidth]{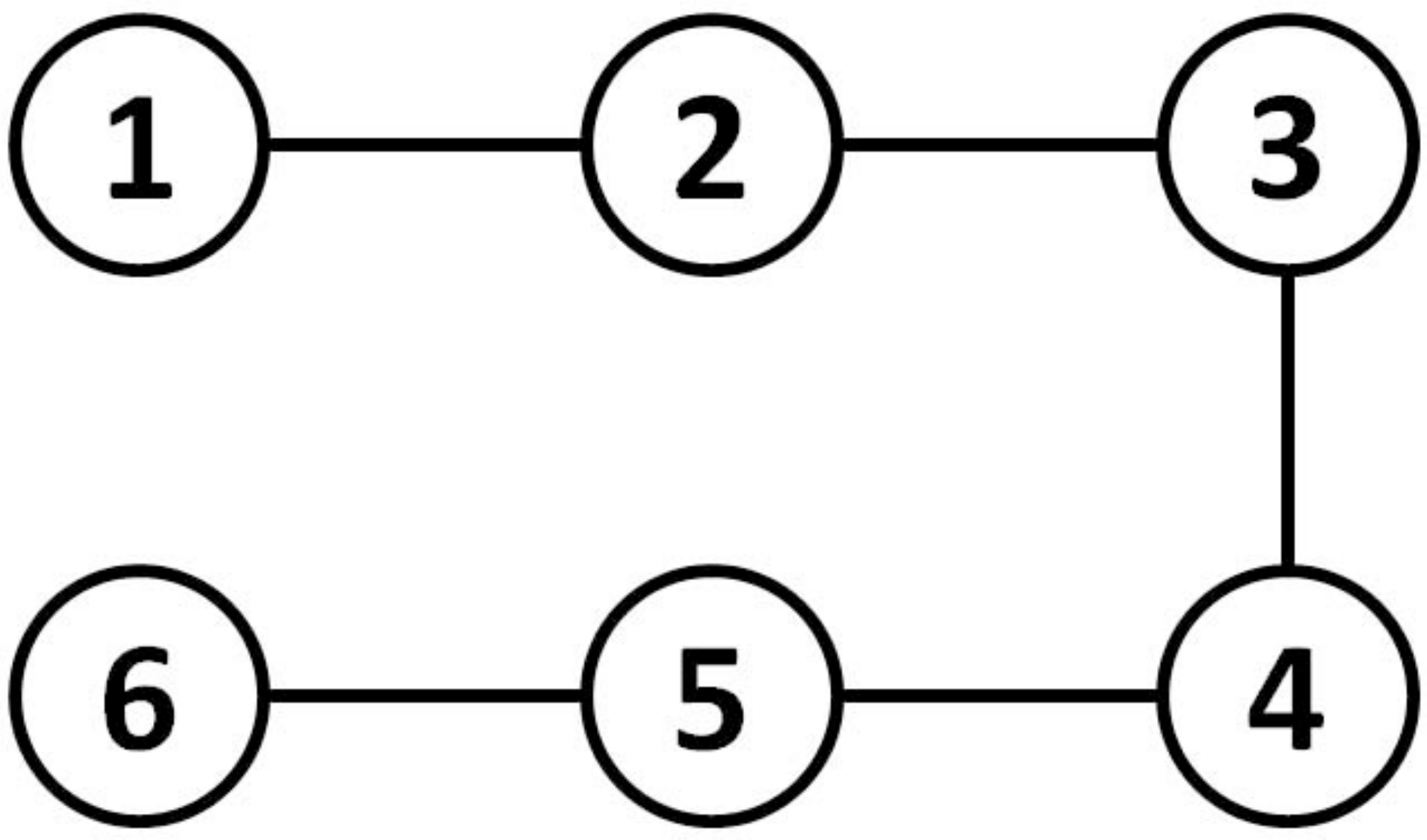}}
\caption{Examples of network topologies}
\end{figure}

For example as shown in Figure \ref{fig:UndirectGraph1}, where a graph of 6 nodes is shown.
The neighborhood of node $5$ is $N_5=\{1,3,4,6\}.$
As node $5$ has the most number of neighbors, the maximum degree $\Delta_{\textrm{max}}=\Delta_5=4$.
Another graph which also has the maximum degree $\Delta_{\textrm{max}}=4$ is shown in Figure \ref{fig:UndirectGraph2}.
A possible spanning tree of the graph shown in Figure \ref{fig:UndirectGraph1} is shown in Figure \ref{fig:SpanningTree}, in which the maximum degree of the graph is 2.
Note that a small maximum degree implies that the graph is sparse.

\subsection{System Dynamics and Sensing Model}

We consider the linear dynamic system and sensing model, as
\begin{equation}\label{Equ:SysDynamics}
\begin{split}
x[k]&=A[k] x[k-1] + B[k] w[k],\\
z[k]& = H[k] x[k] + v[k],
\end{split}
\end{equation}
where $x[k]\in\mathbb{R}^n$ and $z[k]\in\mathbb{R}^m$ are, respectively, the state and measurement vector.
$A[k]\in\mathbb{R}^{n\times n}$ is the state transition matrix.
$B[k]\in\mathbb{R}^{n\times p}$ models the process noise.
$H[k]\in\mathbb{R}^{m\times n}$ is the state-to-measurement matrix.
$w[k]\in\mathbb{R}^{p}$ and $v[k]\in\mathbb{R}^{m}$ are, respectively, the process and measurement noise at time instant $k$.
Both $w[k]$ and $v[k]$ are assumed to be white Gaussian noise, where $w[k]\sim\mathcal {N}(0,Q[k])$ and $v[k]\sim\mathcal {N}(0,R[k])$.

Through this paper, for the purpose of simplicity, we drop the time instant $k$.
As we only care about two relative time instants, namely, the current and the next time instant, we use none sub/superscript and superscript ``$+$", to denote them, respectively.
For example, if the current time instant is $k-1$, the first equation in \eqref{Equ:SysDynamics} can be written as $x^+=A^+x+B^+w^+$.

\subsection{Weighted Least Square Estimator and Kalman Filter}\label{Sec:WLStoKF}
We briefly mention the derivation of the Kalman filter from the recursive Weighted Least Square (WLS) estimator\cite{farrell2008aided}.
Given a measurement at a certain time instant \[z=Hx+v,\] the WLS estimate $\hat{x}_\text{wls}$ is
\begin{equation}\label{Equ:WLS_Estimator}
\hat{x}_\text{wls}=P_\text{wls}H^\top Wz,
\end{equation}
where $P_\text{wls} = (H^\top WH)^{-1}$ is the estimation error covariance and $W=R^{-1}$ is the weighting matrix.

Suppose that prior to obtaining the measurement $z$, the prior estimate $\bar{x}$ is known as $\mathbb{E}\langle\bar{x}\rangle=x$ and $\textrm{var}(\bar{x})=P$.
Let $\bar{\eta}$ be the prior estimation error defined as $\bar{\eta} \triangleq \bar{x}-x$.
It follows that $\mathbb{E}\langle\bar{\eta}\bar{\eta}^\top\rangle=P$.
By organizing all the available information for estimating $x$ together, a hypothetical measurement model can be written as
\[\tilde{z}=\tilde{H}x+\tilde{v},\]
where $\tilde{z}\triangleq[z^\top ,\bar{x}^\top]^\top$, $\tilde{v}\triangleq[v^\top ,\bar{\eta}^\top]^\top$ and $\tilde{H}\triangleq[H^\top ,\Id_n^\top]^\top$ with $\Id_n\in\mathbb{R}^{n\times n}$ being the identity matrix.
Assume that $\mathbb{E}\langle v\bar{\eta}^\top\rangle=0$.
The measurement noise covariance for the above hypothetical model, defined as $\tilde{R}\triangleq\mathbb{E}\langle \tilde{v}\tilde{v}^\top\rangle$, can be written as $\tilde{R}=\textrm{blkdiag}(R,P)$, where $\textrm{blkdiag}(\cdot)$ denotes the block diagonal matrix.
Let $\tilde{W}=\tilde{R}^{-1}$.
By applying \eqref{Equ:WLS_Estimator} with $\tilde{H}$, $\tilde{z}$ and $\tilde{W}$ playing, respectively, the role of $H$, $z$, and $W$, the posterior state estimate and covariance can be obtained, respectively, as
\begin{equation}\label{Equ:Kalman:Update}
\begin{split}
\hat{x}&=M\tilde{H}^\top \tilde{R}^{-1}\tilde{z}
=M(H^\top R^{-1}z + P^{-1}\bar{x}),\\
M &= \left(\tilde{H}^\top\tilde{R}^{-1}\tilde{H}\right)^{-1}
=(H^\top R^{-1}H+P^{-1})^{-1}.
\end{split}
\end{equation}
Eq. \eqref{Equ:Kalman:Update} is commonly referred as the update steps of the Kalman filter \cite{Kalman1960}.
The prediction steps of the Kalman filter are
\begin{equation}\label{Equ:Kalman:Prediction}
\begin{split}
x^+ &=A\hat{x},\\
P^+ &=AMA^\top+BQB^\top.
\end{split}
\end{equation}
By recursively applying \eqref{Equ:Kalman:Update} and \eqref{Equ:Kalman:Prediction} at each time instant, the Kalman filter is able to estimate the state of the linear dynamic system \eqref{Equ:SysDynamics}, in an optimal manner.
It is well known that this is guaranteed to be achieved if and only if $(A,H)$ is an observable pair for all time instants.

\section{Distributed Implementation}

Suppose that there are $N$ nodes.
We assume that node $i$ is able to obtain the local measurement $z_i\in\mathbb{R}^{m_i}$ such that \[z_i=H_ix+v_i,\] where $H_i$ is the local state-to-measurement matrix and can be different for each node, and $v_i\sim\mathcal{N}(0,R_i)$ is the local measurement noise.
Let $\bar{x}_i $ and $\hat{x}_i$ be the local prior and posterior estimates of $x$, respectively.
Let $\bar{\eta}_i\triangleq\bar{x}_i -x$ and $\hat{\eta}_i\triangleq\hat{x}_i-x$ be the local prior and posterior estimation error, respectively.
Let $P_i\triangleq\mathbb{E}\langle \bar{\eta}_i \bar{\eta}_i^\top \rangle$ and $M_i\triangleq\mathbb{E}\langle\hat{\eta}_i\hat{\eta}_i^\top\rangle$  be the local prior and posterior estimation error covariance, respectively.
By using only its own measurement at every time instant, each node can recursively implement a local WLS estimator.
By using its own measurement and local prior estimate, each node can derive the local Kalman filter from the local WLS estimator in the same manner as shown in Section \ref{Sec:WLStoKF} as
\begin{equation}\label{Equ:LKF}
\begin{split}
\hat{x}_i&=M_i(H_i^\top R_i^{-1}z_i + P_i^{-1}\bar{x_i}),\\
M_i &= (H_i^\top R_i^{-1}H_i+P_i^{-1})^{-1},\\
x^+_i &=A\hat{x_i},\\
P^+_i &=AM_iA^\top+BQ_iB^\top.
\end{split}
\end{equation}
Unfortunately, such a filter does not make use of the information flow in the network.

In this section, we start with a novel local WLS estimator, where each node uses information from its own measurement and prior estimate, and information from its local neighbors' measurements and prior estimates obtained through local communications.
With some proper approximations in the derivation of the covariance matrix, we then derive the update steps of our proposed distributed implementation.

\subsection{Local WLS Estimator with Generalized Measurement}
Assume that at each time instant, each node, say node $i$, is able to get access to its local neighbors' prior local estimates and estimation error covariances, as well as their local measurements, and local measurement noise covariance, i.e., $\bar{x}_j$, $P_j$, $z_j$ and $R_j$, $\forall j\in N_i$, through local communication.
We let each node form its local \emph{generalized measurement}, denoted as $\tilde{z}_i$, using not only its own measurement and prior estimate of the state, but also its local neighbors'.
That is, let $\tilde{z}_i$ be defined as \[\tilde{z}_i\triangleq[z_i^\top,z_{j_1}^\top,\ldots,z_{j_{|N_i|}}^\top,\bar{x}_{i}^\top,\bar{x}_{j_1}^\top,\ldots,\bar{x}_{j_{|N_i|}}^\top]^\top,\] where $ j_\ell\in N_i$.
Note that $\bar{x}_i=x+\bar{\eta}_i$.
Let $\mathcal{H}_i\in\mathbb{R}^{(n|J_i|)\times n}$ be the matrix formed by stacking $|J_i|$ $n$ by $n$ identity matrices.
Let $\tilde{\bar{\eta}}_i\triangleq[\bar{\eta}_i^\top,\bar{\eta}_{j_1}^\top,\ldots,\bar{\eta}_{j_{|N_i|}}^\top ]^\top$.
The sensing model of node $i$ with respect to is local generalized measurement can be written as \[\tilde{z}_i=\tilde{H}_ix+\tilde{v}_i,\]
where $\tilde{H}_i=[H_i^\top,H_{j_1}^\top,\ldots,H_{j_{|N_i|} }^\top,\mathcal{H}_i^\top]^\top$, and $\tilde{v}_i=[v_i^\top,v_{j_1}^\top,\ldots,v_{j_{|N_i|} }^\top,\tilde{\bar{\eta}}_i^\top]^\top$.
Assume that $\mathbb{E}\langle v_i\bar{\eta}_j^\top\rangle=0,\ \forall i,j$.
It follows that the generalized measurement noise covariance for node $i$, defined by $\tilde{R}_i\triangleq\mathbb{E}\langle\tilde{v}_i\tilde{v}_i^\top\rangle$, can be written as $\tilde{R}_i=\textrm{blkdiag}(R_i,R_{j_1},\ldots,R_{j_{|N_i|}},\mathcal{P}_i)$, where
\begin{equation}\label{Equ:P_i}
\Pcal_i\triangleq\mathbb{E}\langle\tilde{\bar{\eta}}_i\tilde{\bar{\eta}}_i^\top\rangle=\left[
\begin{array}{cccc}
  P_{i}&P_{ij_1}&\cdots&P_{i\jN}\\
 P_{j_1i}&\ddots&&\vdots\\
 \vdots&&\ddots&\vdots\\
 P_{\jN i}&\cdots&\cdots&P_{\jN}
\end{array}\right].
\end{equation}
In \eqref{Equ:P_i}, 
$P_{ij}\triangleq\mathbb{E}\langle\bar{\eta}_i\bar{\eta}_j^\top\rangle$, where $i\neq j$.
By substituting the local generalized variables $\tilde{z}_i$, $\tilde{H}_i$ and $\tilde{R}_i$ into \eqref{Equ:Kalman:Update}, each node is able to update its local posterior estimate and estimation error covariance as
\begin{equation}\label{Equ:LocalKalman:Update:NoApproximation}
\begin{split}
\hat{x}_i&=M_i\left(
\sum_{j\in J_i}{H_j^\top R_j^{-1}z_j} + \Hcal^\top \Pcal_i^{-1}\tilde{\bar{x}}_i
\right),\\
M_i &=\left(\sum_{j\in J_i} H_j^\top R_j^{-1} H_j +\Hcal^\top \Pcal_i^{-1}\Hcal\right)^{-1},
\end{split}
\end{equation}
where $\tilde{\bar{x}}_i=[\bar{x}_{i}^\top,\bar{x}_{j_1}^\top,\ldots,\bar{x}_{j_{|N_i|}}^\top]^\top$.

\subsection{Approximation of $\Pcal_i^{-1}$}
Due to the structure of $\Pcal_i$ as shown in \eqref{Equ:P_i}, the computation of $\Pcal_i^{-1}$ requires the knowledge of cross-correlation $P_{ij_\ell}$, which is nonzero in general, $\forall j_\ell\in N_i$.
As studied in \cite{olfati2009kalman,kamal2013ICFTAC}, the computation of each $P_{ij_\ell}, j_\ell\in N_i$ will require the knowledge of node $j_\ell$'s neighbors, and so on so forth.
Therefore, $P_{ij_\ell}$ cannot be computed locally.
This makes the update steps in \eqref{Equ:LocalKalman:Update:NoApproximation} not distributed.
Therefore, an appropriate approximation to compute $\Pcal_i$ in a distributed manner is required.
A common way is to assume $\mathbb{E}\langle\bar{\eta}_i\bar{\eta}_j^\top\rangle=0,\forall i\neq j$.
In such a case, $\Pcal_i$ will become a block diagonal matrix and hence $\Pcal_i^{-1}\approx\textrm{blkdiag}(P_i^{-1},P_{j_i}^{-1},\ldots,P_{\jN}^{-1})$.
Therefore, \eqref{Equ:LocalKalman:Update:NoApproximation} can be approximated as
\begin{equation}\label{Equ:LocalKalman:Update:Approximation1}
\begin{split}
\hat{x}_i&=M_i\left(
\sum_{j\in J_i}{H_j^\top R_j^{-1}z_j} + \sum_{j\in J_i}P_j^{-1}{\bar{x}}_j
\right),\\
M_i &=\left(\sum_{j\in J_i} H_j^\top R_j^{-1} H_j +  \sum_{j\in J_i}P_j^{-1}\right)^{-1}.
\end{split}
\end{equation}

Compared with the local Kalman filter \eqref{Equ:LKF} that uses only information from its own measurement and prior estimate, \eqref{Equ:LocalKalman:Update:Approximation1} uses information from both its own and its neighbors' measurements and prior estimates.
Specifically,
$P_i^{-1}\bar{x}_i$ (respectively, $P_i^{-1}$) in \eqref{Equ:LKF} is substituted by $\sum_{j\in J_i}P_j^{-1}\bar{x}_j$ (respectively, $\sum_{j\in J_i}P_j^{-1}$) in \eqref{Equ:LocalKalman:Update:Approximation1}.
As each node iterates enough, ideally, it would have a good understanding about the target.
Then the local estimation error covariance should be close to the hypothetical centralized estimation error covariance, denoted as $P_{\text{c}}$, i.e., $P_i\rightarrow P_{\text{c}}$.
In such a case, $\sum_{j\in J_i}P_j^{-1}\rightarrow|J_i|P_{\text{c}}^{-1}$.
That is, the local information matrices (i.e., inverse of the local estimation error covariance matrices) from all nodes in the inclusive neighborhood $J_i$, are overused by a factor of $|J_i|$.
In other words, the uncertainties from all nodes in the inclusive neighborhood are underestimated by a factor of $|J_i|$.
To avoid such underestimations, we modify the approximated $\Pcal_i$ as $\Pcal_i\approx |J_i|\textrm{blkdiag}(P_i,P_{j_1},\ldots,P_{\jN})$.
With this modified approximation of $\Pcal_i$, \eqref{Equ:LocalKalman:Update:Approximation1} can be revised as
\begin{equation}\label{Equ:LocalKalman:Update:Approximation2}
\begin{split}
\hat{x}_i&=M_i\left(
\sum_{j\in J_i}{H_j^\top R_j^{-1}z_j} + \frac{1}{|J_i|}\sum_{j\in J_i}P_j^{-1}{\bar{x}}_j
\right),\\
M_i &=\left(\sum_{j\in J_i} H_j^\top R_j^{-1} H_j +  \frac{1}{|J_i|}\sum_{j\in J_i}P_j^{-1}\right)^{-1}.
\end{split}
\end{equation}
We refer \eqref{Equ:LocalKalman:Update:Approximation2} as the update steps of the Information-driven Fully Distributed Kalman filter (IFDKF).
The prediction steps are
\begin{eqnarray*}
\nonumber
P_i^+ &=& A M_i {A}^\top + B Q{B }^\top,\\
\nonumber
\bar{x}_i^+ &=& A\hat{x}_i.
\end{eqnarray*}
\section{Information-driven Fully Distributed Kalman Filter}

The details of the proposed IFDKF are summarized in Algorithm \ref{Algorithm:IFDKF}.
\begin{algorithm}\label{Algorithm:IFDKF}
\DontPrintSemicolon 
\KwIn{Initial guess of state $x_i^0$ and estimation error covariance $P_i^0$}
\KwOut{Local posterior state estimate $\hat{x}_i$ and posterior estimation error covariance $M_i$}
\If{$k=1$}
{$\bar{x}_i=x_i^0,P_i=P_i^0$}
takes the local measurement $z_i$\\
Computes $S_i\gets H_i^\top R_i^{-1}H_i,\quad y_i\gets H_i^\top R_i^{-1}z_i$\\
Sends $S_i$, $y_i$, $\bar{{x}}_i$ and $P_i$ to node $j$, $\forall j\in N_i$\\
Receives $S_j$, $y_j$, $\bar{{x}}_j$ and $P_j$ from node $j$, $\forall j\in N_i$\\
Forms the inclusive neighborhood $J_i\gets N_i\cup\{i\}$\\
Computes $\bar{y}_i\gets\sum_{j\in J_i}{y_i},\quad \bar{S_i}\gets\sum_{j\in J_i}{S_i}$\\

Updates the local state estimate and estimation error covariance
\begin{eqnarray}
\label{Equ:IFDKF:StateUpdate}
\hat{x}_i&\gets& M_i\left(
\bar{y}_i + \frac{1}{|J_i|}\sum_{j\in J_i}P_j^{-1}{\bar{x}}_j\right)\\
\label{Equ:IFDKF:CovUpdate}
M_i &\gets&\left(\bar{S}_i +  \frac{1}{|J_i|}\sum_{j\in J_i}P_j^{-1}\right)^{-1}
\end{eqnarray}
\label{Alg:Update}
\Return{$\hat{x}_i,M_i$}\\
Predicts the local prior estimate and estimation error covariance for next time instant
\begin{eqnarray}
\nonumber
\bar{x}_i^+&\gets& A\hat{x}_i\\
\nonumber
P_i^+ &\gets& A M_i {A}^\top + B Q{B }^\top
\end{eqnarray}
\caption{IFDKF at Node $i$ at Time Instant $k$}
\end{algorithm}
The advantages of this algorithm will be discussed in detail in this section.
\subsection{Robustness against the Presence of Naivety}

\begin{definition}[Naive Node]
Let $\mathfrak{N}$ be the set of indices whose corresponding nodes are \emph{naive}.
Then $i\in\mathfrak{N}$ if $(A,H_{J_i})$ is not an observable pair, where $H_{J_i}$ is the matrix obtained by stacking all $H_j$ such that $j\in J_i$.
\end{definition}

As each node obtains its local measurement, communicates with its neighbors and updates its local estimate at the same frequency, the naive nodes can neither directly nor indirectly obtain any complete information related to the current state.
They can only be indirectly influenced by some time-delayed complete information related to the past state.
Note that \eqref{Equ:IFDKF:StateUpdate} can be further written as
\begin{equation}\label{Equ:Derivation1}
\hat{x}_i= \bar{x}_i + M_i\left(-M_i^{-1} \bar{x}_i
+\bar{y}_i + \frac{1}{|J_i|}\sum_{j\in J_i}P_j^{-1}{\bar{x}}_j\right).
\end{equation}
From \eqref{Equ:IFDKF:CovUpdate}, one can obtain that $M_i ^{-1} =\bar{S}_i +  \frac{1}{|J_i|}\sum_{j\in J_i}P_j^{-1}$.
By substituting $M_i ^{-1}$ into \eqref{Equ:Derivation1}, it follows that
\begin{eqnarray}\label{Equ:Update:ConsensusForm}
\hat{x}_i=
\bar{x}_i
+ M_i\left(\bar{y}_i -\bar{S}_i\bar{x}_i\right)
+\frac{1}{|J_i|}M_i\sum_{j\in J_i} P_j^{-1} \left({\bar{x}}_j-\bar{x}_i\right).
\end{eqnarray}
Therefore, \eqref{Equ:IFDKF:StateUpdate} can be equivalently represented as \eqref{Equ:Update:ConsensusForm}.

In \eqref{Equ:Update:ConsensusForm}, both of the last two terms are left-multiplied by $M_i$, whose 2-norm implies the uncertainty of the local posterior estimate of node $i$.
If node $i$ is relatively confident about its own local estimation, it prefers to do less correction compared to its own prior estimate since $\|M_i\|_2$ is relative small and consequently attenuates the effects of the last two terms.
Note that $\bar{y}_i-\bar{S}_i\bar{x}_i=\sum_{j\in J_i}{H_j^\top R_j^\top(z_j-H_j\bar{x}_i)}$.
Therefore if node $i$ is less confident about its own estimation, it will depend more on the latest measurements from its inclusive neighborhood ($z_j,\ \forall j\in J_i$), as implied by the second term of \eqref{Equ:Update:ConsensusForm}.
Meanwhile, it also depends more on its neighbor's prior estimates ($\bar{x}_j,\ \forall j\in N_i$), as implied by the last term of \eqref{Equ:Update:ConsensusForm}.
Moreover, instead of equally weighing the difference between its own and neighbor's prior estimates as in the KCF, each node weighs the difference by the matrix $P_j^{-1}$, for each $j\in N_i$.
The matrix $P_j^{-1}$, known as the prior information matrix of node $j$, will have relative small 2-norm if $\bar{x}_j$ is relatively erroneous.
Therefore, node $i$'s local estimate is pushed more toward a neighbor with more accurate estimate.
For a naive node, its prior estimates will be very erroneous, which implies that the 2-norm of its information matrix will be very small.
Therefore, for each node, if it is not very confident about it own estimate, the last term in \eqref{Equ:Update:ConsensusForm} will wisely drive its local estimate toward its neighbors' prior estimates based on their relative confidence.
In \eqref{Equ:IFDKF:CovUpdate}, the term $\frac{1}{|J_i|} \sum_{j \in J_i}P_j^{-1}$ reflects the fact that due to the communication of local prior estimates between local neighbors, each node's posterior estimation error covariance is also affected by its neighbors' prior estimation error covariances.

If we remove the last term in \eqref{Equ:Update:ConsensusForm} and replace the term $\frac{1}{|J_i|} \sum_{j \in J_i}P_j^{-1}$ with $P_i^{-1}$, we essentially obtain a local Kalman filter that uses information from its own measurement and prior estimate as well as its neighbors' measurements.
However, without the coupling between neighbors' prior estimates in the posterior state estimate update, or the incorporation of neighbors' prior estimation error covariances in the posterior estimation error covariance update, such a local Kalman filter does not work in the presence of naive nodes as the KCF.
Superficially, \eqref{Equ:Update:ConsensusForm} and \eqref{Equ:IFDKF:CovUpdate} might look similar to the KCF.
However, the specific form of the last term in \eqref{Equ:Update:ConsensusForm} and the second term in \eqref{Equ:IFDKF:CovUpdate} play an important role in dealing with the presence of naive nodes.
It is also worth mentioning that the last term in \eqref{Equ:Update:ConsensusForm} and the second term in \eqref{Equ:IFDKF:CovUpdate} are not added or introduced in an add-hoc manner based on the intuition but rather are derived based on the novel local WLS estimator which forms its local generalized measurement by collecting both its own and local neighbors' measurements and their prior estimates.
In contrast, the KCF is derived based on a local Kalman filter that uses the measurements obtained by nodes in its inclusive neighborhood and then a consensus term is added to the posterior state estimate of the local Kalman filter by following the intuition to drive the nodes' local estimates to be consistent.
Note that both the above mentioned local Kalman filter and the KCF rely on the assumption that the target is jointly observed in the inclusive neighborhood of each node (that is, the absence of naive nodes). Otherwise, the naive nodes' estimation error covariances and estimates would diverge and in turn deteriorate those nodes with good estimates.
\subsection{Fully Distributed}

In the GKCF and ICF, the design of the consensus step size $\epsilon$ is required to be chosen between $0$ and $1/\Delta_\text{max}$, where $\Delta_\text{max}$ is the maximum degree of the network graph.
If the maximum degree is changing with time, either due to the switches of the communication topology, or due to some unknown failures of the communication channels, a proper selection of $\epsilon$ might not be as good as before.
Even if $\Delta_\text{max}$ is known, it is not clear how to select a nice $\epsilon$ in the studies of GKCF and ICF.
In \cite{kamal2011GKCF,kamal2013ICFTAC}, after obtaining a new local measurement, each node is allowed to communicate with its local neighbors for infinite times, before it finally updates its local estimate at the current time instant.
Therefore, in such a case, any selection between $0$ and $1/\Delta_{\text{max}}$ will guarantee that every node's local estimate asymptotically becomes consistent before its next measurement update.
Unfortunately, infinite communication steps between measurement steps are not realistic.
In the realistic case where every node takes a measurement, communicates with its local neighbors, and updates its local estimate at the same frequency, the selection of $\epsilon$ becomes more critical as a large $\epsilon$ might make the filter diverge but a small $\epsilon$ might slow down the rate of convergence of each node's local estimate to its neighbors'.

Besides $\Delta_\text{max}$, the ICF also requires each node to know the total number of nodes $N$ in the network in order to asymptotically approach the centralized solution via infinite communication steps before it updates its local estimate.
Similar to the discussion of $\Delta_{\text{max}}$, $N$ could be changing over time due to reasons such as node failures or inclusion of new nodes in the network.
Without knowing the correct $N$, each node will either overestimate or underestimate the uncertainties of its neighbors' prior estimates, and consequently slow down the convergence rate of its local estimate to the true state value.
Moreover, in the considered case of this paper (i.e., the same frequency of obtaining measurement, communication and estimate update), the ICF is not able to obtain the optimal solution anyways even with the correct knowledge of $N$ as with one iteration the embedded consensus algorithm is only able to return the local average instead of the global average.

The proposed IFDKF, different from the consensus-embedded filters, does not require any global information and therefore does not require any parameter determination or tuning.
By using only local information, it is run in an automated manner and is adaptive to the locally unknown changes in the network.

\section{Experimental Evaluation}

In this section, we evaluate the performance of the proposed IFDKF by comparing it with some existing algorithms, namely, the KCF, the GKCF and the ICF.
We adopt the simulation model from \cite{kamal2013ICFTAC}, where a target is moving with a certain model and randomly affected by some stochastic process noise.
Assume that the target is moving in a certain area that is monitored by a network with 6 cameras.
Assume that there exists at least one naive node.
We let the system dynamics and sensing model be linear time invariant as shown in \eqref{Equ:SysDynamics}, where the state vector $x\in\mathbb{R}^4$ contains two position and two velocity components.
Let $B=\Id_4$, $Q=\textrm{diag}(10,10,1,1)$, $R=100\Id_2$, and
\begin{equation*}
  A=\left[
  \begin{array}{cccc}
1&0&1&0\\
0&1&0&1\\
0&0&1&0\\
0&0&0&1
\end{array}
  \right].\
\text{Let }  H_{i}=\left[  \begin{array}{cccc}
1&0&0&0\\
0&1&0&0
\end{array}
  \right]
\end{equation*}
if node $i$ is able to directly observe the target, and $H_{i}=\zeros$ otherwise.
It is obvious that if $i\in\mathfrak{N}$, $H_{J_i}=\zeros$ and therefore $(A,H_{J_i})$ is not an observable pair.
\subsection{Performance in Presence of Naive Nodes}
\begin{figure}
\includegraphics[width=\linewidth]{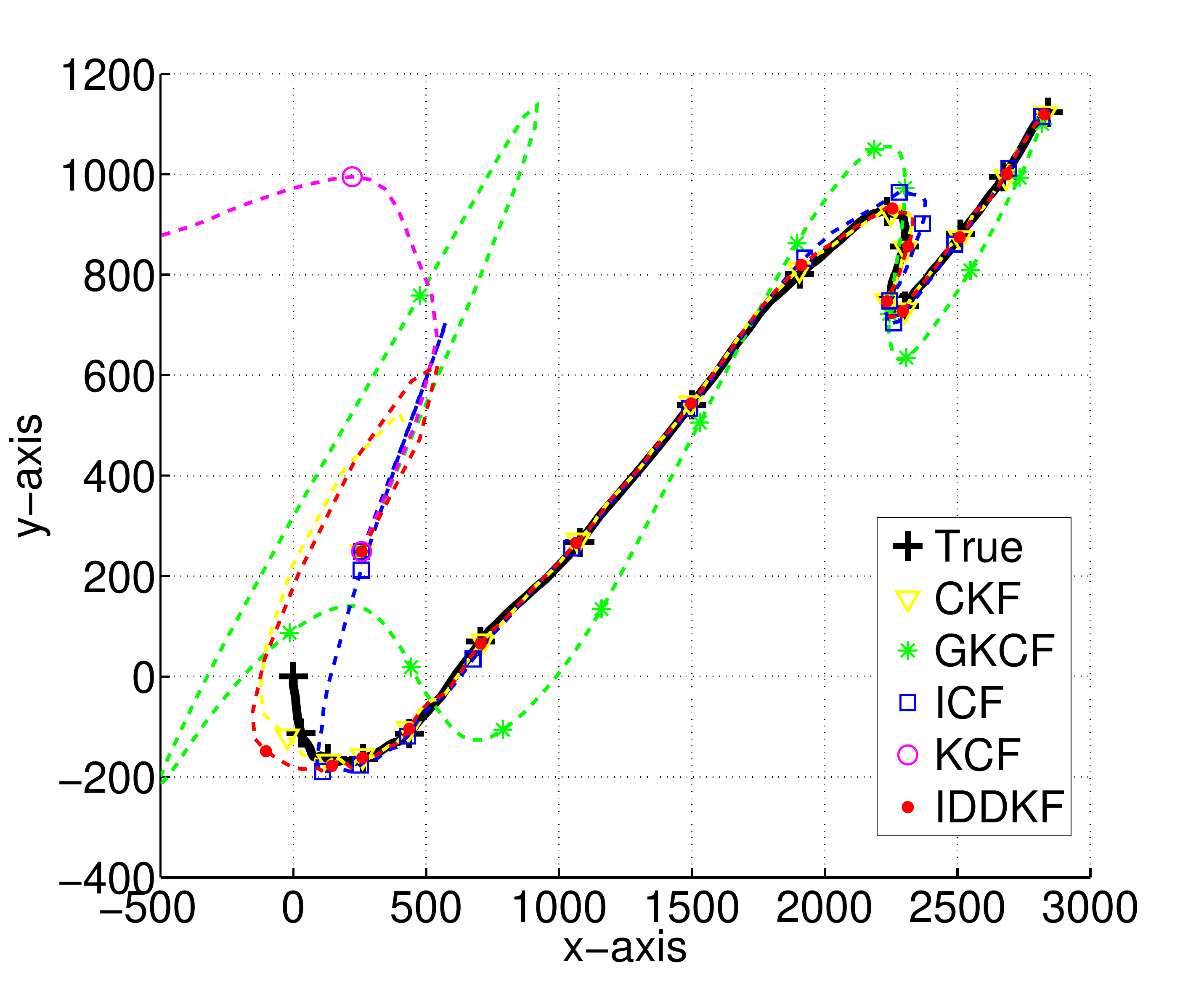}
\caption{Comparison of algorithms: tracking task of a moving target\label{fig:TrackingCompareDenseTopo}}
\end{figure}
We first adopt the graph as shown in Figure \ref{fig:UndirectGraph1}, in which the communication topology is relatively dense.
We only let the target be directly observed by node 1.
Therefore, the naive node index set is $\mathfrak{N}=\{3,4,6\}$.
A tracking task is tested to compare the performance of the KCF, GKCF, ICF and IFDKF.
For each distributed estimation algorithm, the means of the local estimated locations for all 150 time instants are plotted in Figure \ref{fig:TrackingCompareDenseTopo} in different colors.
For each algorithm, a different type of mark is also plotted every $10$ time instants.
The true location of the target and the estimates of the centralized Kalman filter (CKF) are also plotted as references.
The target is initially located at the origin of the coordinate.
The initial local prior estimate of each node is randomly chosen from 0 to 500 with equal probability.
To be fair, the initial local prior estimate of each node is chosen to be consistent for different algorithms.
The initial prior estimate of the CKF is chosen as the mean of the selected initial local prior estimates of all six nodes.
For the algorithms that involve the consensus parameter $\epsilon$, the parameter is selected as $0.65/\Delta_{\textrm{max}}$.
This is consistent with that of the selection in \cite{kamal2013ICFTAC}.

As shown in Figure \ref{fig:TrackingCompareDenseTopo}, the mean of the local estimates obtained by the KCF deviates from the true state quickly.
This is due to the presence of naive nodes in the network.
The detailed analysis of why the KCF is not able to deal with the naive nodes can be found in \cite{kamal2011GKCF}.
Later on in this section, we only focus on the algorithms that are able to tolerate the existence of naive nodes.
It can be observed from Figure \ref{fig:TrackingCompareDenseTopo} that the GKCF is only able to ``roughly" follow the true state.
However, the ICF and the IFDKF perform much better in the sense that they both get close to the centralized solution after about 20 time instants.
It is worth noticing that when there is relatively obvious process noise (the target suddenly changes its direction), the IFDKF reacts more quickly than the ICF.

\subsection{Severe Case: ``Chain" Topology}

In this subsection, we adopt the graph as shown in Figure \ref{fig:SpanningTree}.
We assume that the target is only observed by node $1$.
It follows that $\mathfrak{N}=\{3,4,5,6\}$.
Note that this is a very bad case as it takes 5 iterations for node 6 to be influenced by node 1 indirectly.
Therefore, there will be a delay on the local estimate of node 6.
\begin{figure}
\includegraphics[width=\linewidth]{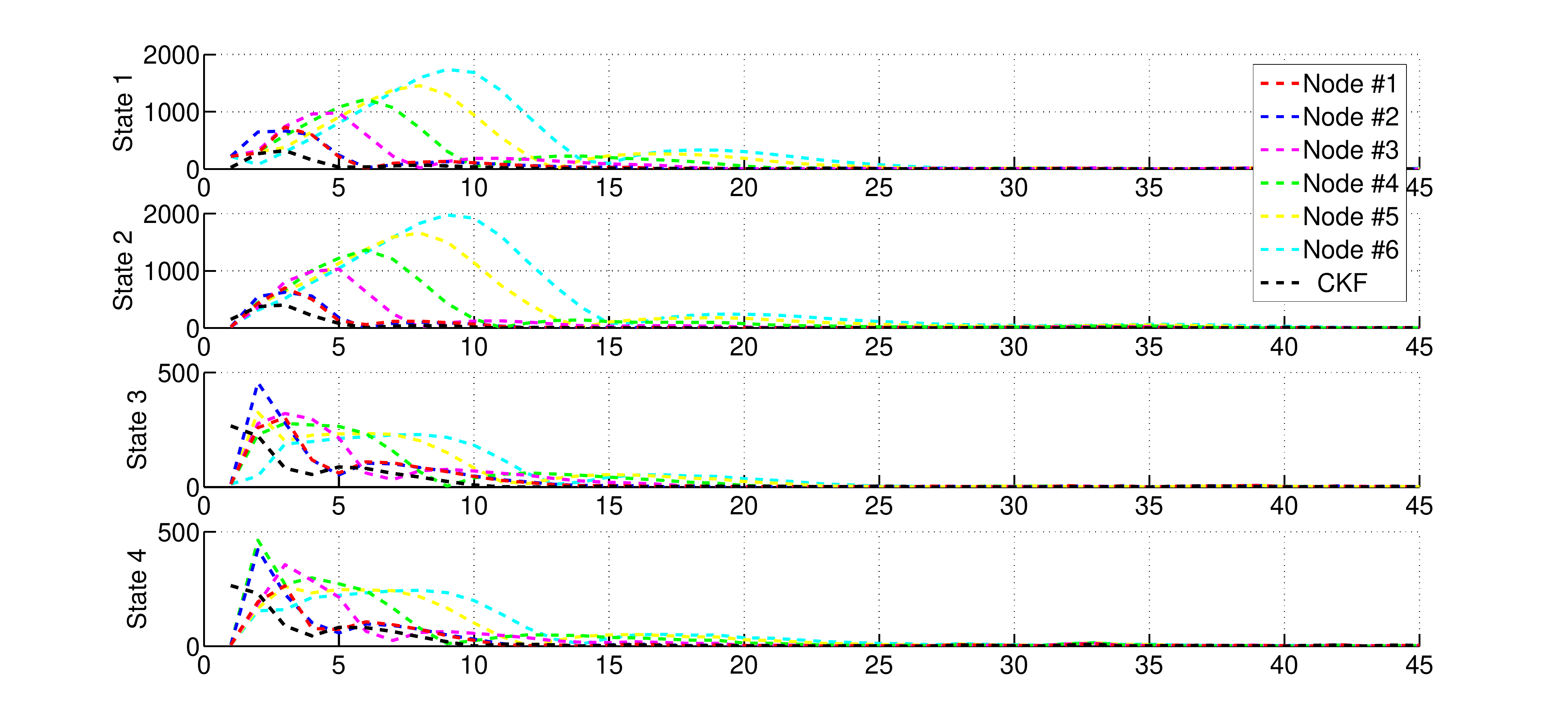}
\caption{$|\hat{\eta}_i|$ of IFDKF for all nodes under ``chain" topology\label{fig:AbsErrIFDKF}}
\end{figure}
The absolute local posterior estimation error of all states at all nodes are shown in Figure \ref{fig:AbsErrIFDKF}.
It is obvious that node 6, which has the longest distance from the node that directly observes the target, has the slowest convergence rate.
Whenever a relative obvious process noise occurs, node 6 is the last one to react.
This phenomenon will significantly become less obvious with the increase of the graph density or the decrease of the maximum distance among all distances from each naive node to its closest non-naive node.

\begin{figure}
\includegraphics[width=\linewidth]{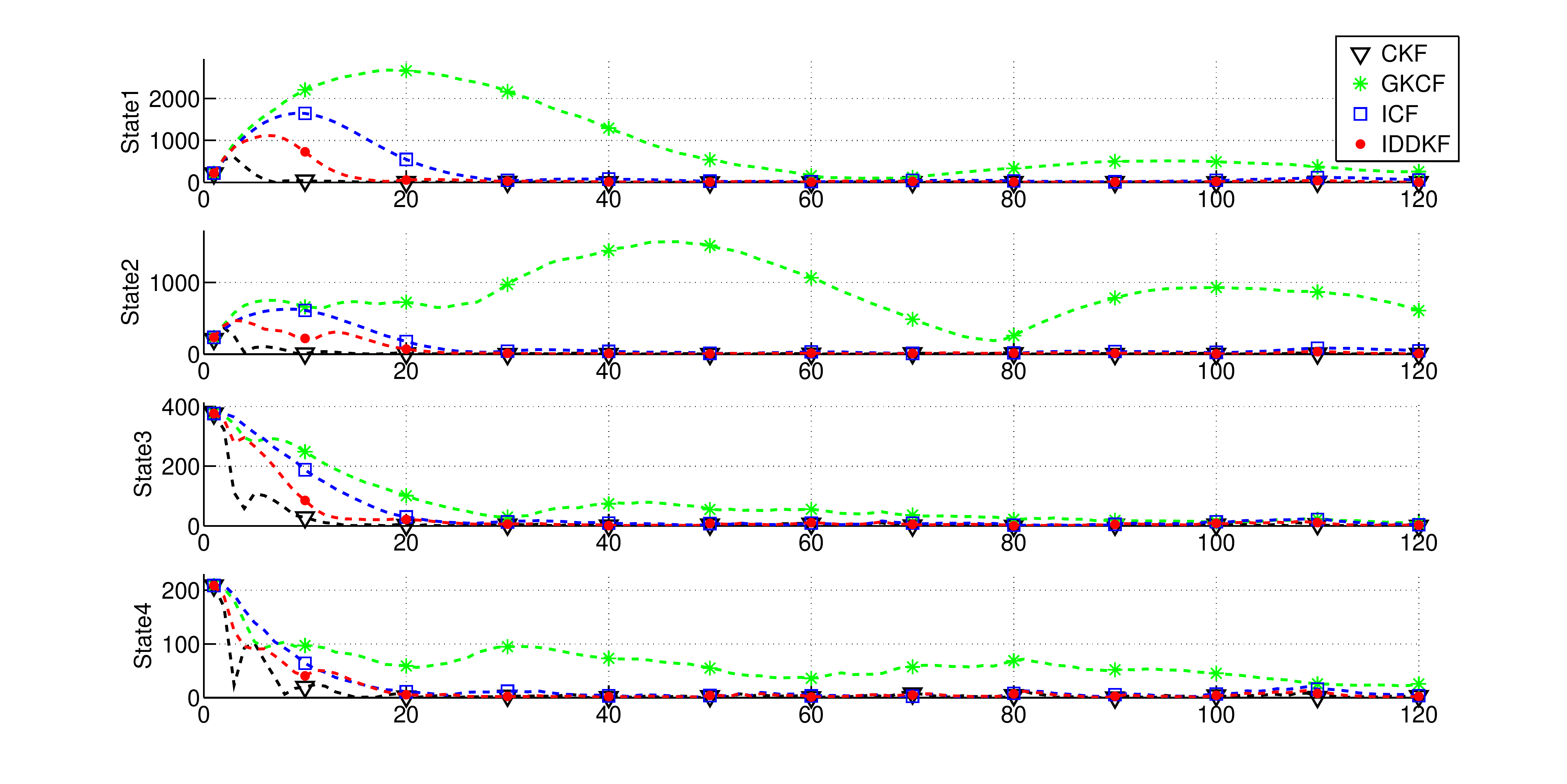}
\caption{Comparison of $\frac{1}{N}\sum_{i=1,\ldots,N}|\hat{\eta}_i|$ from different algorithms\label{fig:AvgErrCompare1}}
\end{figure}

\begin{figure}
\includegraphics[width=\linewidth]{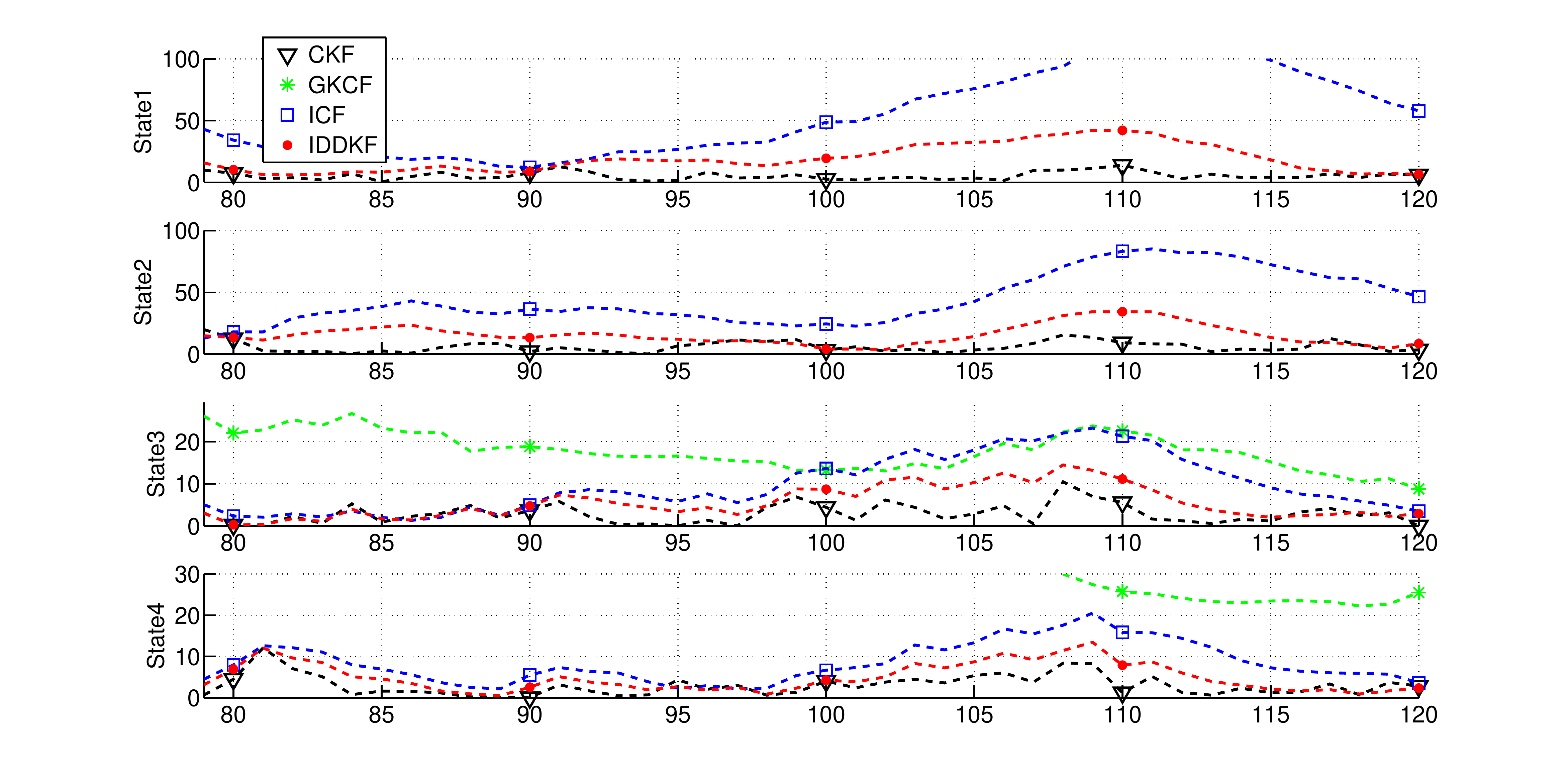}
\caption{Comparison of $\frac{1}{N}\sum_{i=1,\ldots,N}|\hat{\eta}_i|$ in steady state\label{fig:AvgErrCompare2}}
\end{figure}

The mean of the absolute local estimation error for each state among all nodes, i.e., $\frac{1}{N}\sum_{i=1,\ldots,N}|\hat{\eta}_i|$, obtained by the GKCF, ICF and IFDKF are plotted in Figure \ref{fig:AvgErrCompare1}.
As observed, the mean absolute error obtained by the IFDKF converges to the centralized solution more quickly than that of the ICF and the GKCF.
The zoomed-in version of Figure \ref{fig:AvgErrCompare1} is shown in Figure \ref{fig:AvgErrCompare2}.
It is observed that in steady state, in general, the mean absolute error obtained by the IFDKF is lower than that of the other two algorithms.
The IFDKF reacts quicker than the ICF when the target incurs relatively obvious process noise.
For example, the target incurred a disturbance around $k=90$.
Correspondingly, the mean absolute error obtained by each distributed algorithm increases.
It can be observed that the mean absolute error obtained by the IFDKF has the least deviation from the centralized absolute error, and converges back to the centralized one in the quickest rate.
\subsection{Robustness against Changing Global Parameters}
As discussed, the IFDKF is fully distributed and does not require any global information.
Therefore, it automatically adapts to any possible change of the global parameters, and consistently maintains proper performance.
The ICF was shown to outperform the KCF and the GKCF in \cite{kamal2013ICFTAC}.
It does require, however, the knowledge of the graph maximum degree $\Delta_{\text{max}}$ and the total number of nodes $N$ in the network.
In this subsection, we compare the performance of our algorithm with the ICF, when each of the aforementioned global parameters change during the tracking task.
The first test is done on the network whose topology is switched from the one shown in Figure \ref{fig:UndirectGraph2} to the one shown in Figure \ref{fig:SpanningTree} at $k=65$.
The target is assumed to be only directly observed by node 1 all the time.
Therefore, $\Delta_{\text{max}}=4$ for $k<65$ and $\Delta_{\text{max}}=2$ for $k\geq 65$.
\begin{figure}
\includegraphics[width=\linewidth]{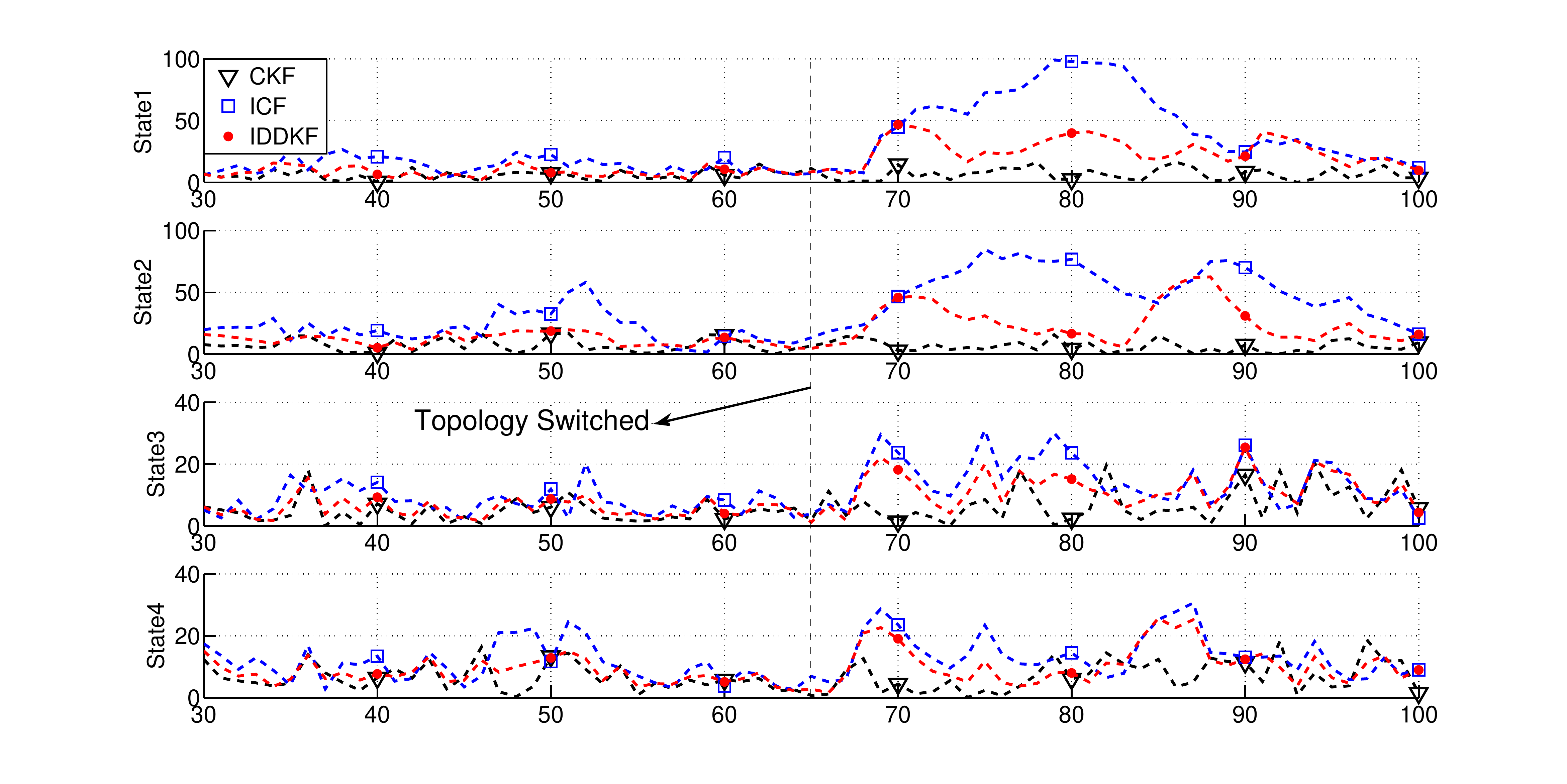}
\caption{Comparison of $\frac{1}{N}\sum_{i=1,\ldots,N}|\hat{\eta}_i|$ with topology changed at $k=65$\label{fig:SwitchingTopology}}
\end{figure}
The comparison of the mean absolute error obtained by the IFDKF and the ICF is shown in Figure \ref{fig:SwitchingTopology}.
As observed the performance of both algorithms perform worse after the topology switches.
This is due to the increased number of naive nodes as $\mathfrak{N}=\{5,6\}$ for $k<65$ and $\mathfrak{N}=\{3,4,5,6\}$ for $k\geq 65$.
Moreover, the difference in performance of the ICF and the IFDKF become more obvious after the switching instant.
This is because the consensus parameter $\epsilon$ selected for the previous topology is no longer appropriate for the current topology.
In this case, since $\Delta_{\text{max}}$ decreases by one half, the selected $\epsilon$ for the graph before the switching instant is only one half of the ``ideal" one after switching.
Therefore, in the ICF, without knowing the change of $\Delta_{\text{max}}$ and keeping using the same $\epsilon$ as before, the naive nodes, whose estimates only depend on the consensus term, will have slower convergence rate and consequently result in worse performance.

The second test is done on the network where some nodes fail during the tracking task.
The original topology is shown in Figure \ref{fig:UndirectGraph2}.
We let the target be directly observed by node 2 and 3.
Therefore $\mathfrak{N}=\{5,6\}$.
At $k=65$, we let node 5 and 6 fail, in the sense that they do not exchange any information with any other node.
\begin{figure}
\includegraphics[width=\linewidth]{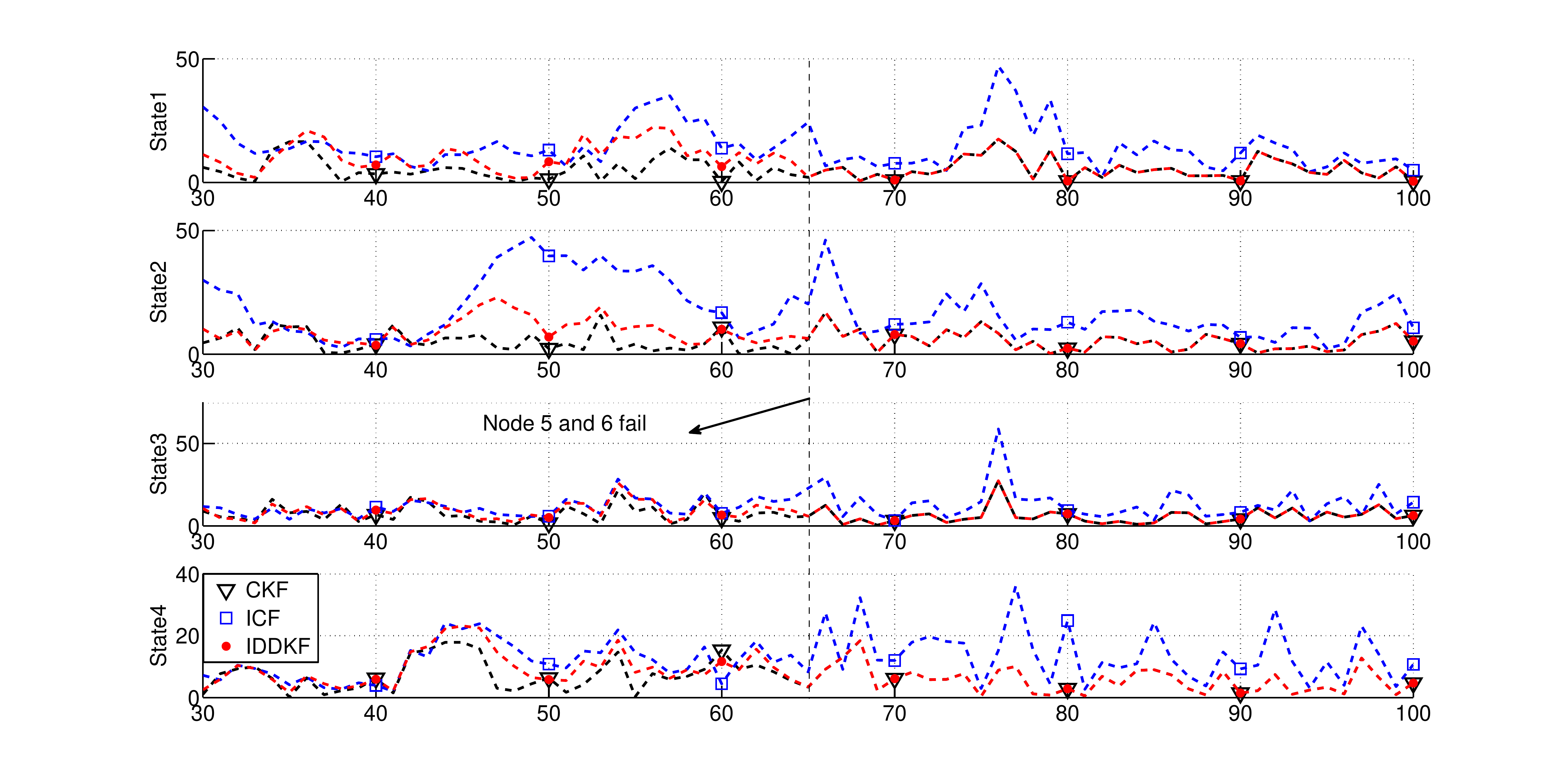}
\caption{Comparison of $\frac{1}{N}\sum_{i=1,\ldots,N}|\hat{\eta}_i|$ with two nodes fail at $k=65$\label{fig:FailedNodes}}
\end{figure}
The comparison of two algorithms in such a case is shown in Figure \ref{fig:FailedNodes}, in which the mean absolute error is computed among all 6 nodes before $k=65$ and only the 4 normal nodes after the failures of node 5 and 6.
Note that when node 5 and 6 fail, there exists no naive node in the graph formed by the rest four nodes as they form an all-to-all topology.
Therefore, the performance of both algorithms, if given the accurate required information, should be identical to the centralized solution as each node is able to obtain as much information as the CKF is able to.
As observed in Figure \ref{fig:FailedNodes}, when $k>65$, the solution obtained by the IFDKF overlaps with the centralized solution.
The ICF, however, fails to converge to the centralized solution.
This is due to the fact that in the ICF, each node, without knowing the decrease on the total number of the nodes in the network, will overestimate its neighbors confidence and therefore slow down the convergence rate.

\section{Conclusion}
We considered the distributed Kalman filtering problem of sensor networks where there exist naive nodes.
We proposed the Information-driven Fully Distributed Kalman filter, in which each node takes measurement, communicates with its local neighbors, and updates its local estimates and estimation error covariance, at the same frequency.
The proposed algorithm does not require any global parameter while is able to handle naive nodes.
The experimental results are used to show the advantages of the proposed algorithm  especially when there exist changes on the global parameters of the network topology.


\addtolength{\textheight}{-12cm}   



\begin{thebibliography}{10}



\bibitem{Kalman1960}
R.~E. Kalman, ``A new approach to linear filtering and prediction problems,''
  \emph{Transactions of the ASME, Journal of Fluids Engineering}, vol.~82,
  no.~1, pp. 35--45, 1960.

\bibitem{RenBeardKingston05}
W.~Ren, R.~W. Beard, and D.~B. Kingston, ``Multi-agent {K}alman consensus with
  relative uncertainty,'' in \emph{Proceedings of the American Control
  Conference}, Portland, OR, June 2005, pp. 1865--1870.

\bibitem{AlighanbariUnbiasedKalman}
M.~Alighanbari and J.~How, ``An unbiased kalman consensus algorithm,'' in
  \emph{Proceedings of the American Control Conference}, Minneapolis, MN, June
  2006, pp. 3519--–3524.

\bibitem{olfati2007distributedKalman}
R.~Olfati-Saber, ``Distributed kalman filtering for sensor networks,'' in
  \emph{Proceedings of the {IEEE} Conference on Decision and Control}, New
  Orleans, LA, December 2007, pp. 5492--5498.

\bibitem{olfati2009kalman}
------, ``Kalman-consensus filter: Optimality, stability, and performance,'' in
  \emph{Proceedings of the {IEEE} Conference on Decision and Control},
  Shanghai, China, December 2009, pp. 7036--7042.

\bibitem{bai2011distributed}
H.~Bai, R.~A. Freeman, and K.~M. Lynch, ``Distributed kalman filtering using
  the internal model average consensus estimator,'' in \emph{Proceedings of the
  American Control Conference}, San Francisco, CA, 2011, pp. 1500--1505.

\bibitem{kamal2011GKCF}
A.~T. Kamal, C.~Ding, B.~Song, J.~A. Farrell, and A.~Roy-Chowdhury, ``A
  generalized kalman consensus filter for wide-area video networks,'' in
  \emph{Proceedings of the IEEE Conference on Decision and Control and European
  Control Conference (CDC-ECC)}, Orlando, FL, 2011, pp. 7863--7869.

\bibitem{kamal2013ICFTAC}
A.~T. Kamal, J.~A. Farrell, and A.~K.~R. Chowdhury, ``Information weighted
  consensus filters and their application in distributed camera networks.''
  \emph{{IEEE} Transactions on Automatic Control}, vol.~58, no.~12, pp.
  3112--3125, 2013.

\bibitem{casbeer2009distributed}
D.~W. Casbeer and R.~Beard, ``Distributed information filtering using consensus
  filters,'' in \emph{Proceedings of the American Control Conference}, St.
  Louis, MO, June 2009, pp. 1882--1887.

\bibitem{farrell2008aided}
J.~A. Farrell, \emph{Aided navigation: GPS with high rate sensors}.\hskip 1em
  plus 0.5em minus 0.4em\relax New York City, NY: McGraw-Hill, 2008.

\end{thebibliography}
\end{document}